\documentclass[a4paper,11pt]{article}
\usepackage[applemac]{inputenc} 
\usepackage{amsfonts,amsmath,amssymb,graphics,amsthm,amscd}
\usepackage[french]{babel}

\usepackage[T1]{fontenc}

\newcommand{\N}{\mathbb{N}}

\newcommand{\Z}{\mathbb{Z}}
\newcommand{\R}{\mathbb{R}}
\newcommand{\Q}{\mathbb{Q}}
\newcommand{\T}{\mathbb{T}}

\def\le{\leqslant}
\def\ge{\geqslant}

\swapnumbers

\theoremstyle{definition}
\newtheorem{defi}{d\'efinition}
\newtheorem{rem}[defi]{- Remarque}

\newtheorem{ex}[defi]{- Exemple}

\theoremstyle{plain}
\newtheorem{prop}{Proposition}
\newtheorem{lem}{Lemme}
\newtheorem{theo}{- Th\'eorème}
\newtheorem{cor}{Corollaire}

\title{L'Algèbre tropicale comme algèbre de la caractéristique  1:\\
 Algèbre linéaire sur les semi-corps idempotents }
\author{Dominique Castella\\
Laboratoire de Mathématiques\\
Université de la Réunion\\
\small{dominique.castella@univ-reunion.fr}}

 \date{17 juillet 2008}
 
\begin{document}
 
\maketitle

\begin{center}
 \end{center}
 
 \begin{abstract}
 We define a formal framework for the study of algebras of type Max-plus, 
Min-Plus, tropical algebras, and more generally algebras over a commutative
idempotent semi-field. This work is motivated by the increasingly diversified 
use of these algebras which occur also in control theory, automata theory  
as well as in algebraic geometry, and in more specific ways in other parts of 
mathematics such as the theory of monoids. 

In this first article, we expecially re-examine linear algebra over idempotent 
semi-fields: the most delicate, but undoubtedly the most interesting point is 
the notion of a singular point seen as a generalization of the notion of 
zero.  We thus rediscover many notions of regularity already introduced for 
matrices, and this permits us to define further notions, new in this context, 
such as that of the kernel of a linear form, and to apply duality to obtain a 
good notion of tropical dimension of a submodule. 

\end{abstract}
\hskip1cm {\small  {\it Keywords} : Max-plus algebra, tropical algebra, idempotent semi-fields.}

 \vskip1cm
 
 \section{Introduction}
 
 Nous définissons un cadre formel pour l'étude des algèbres de types Max-Plus, Min-Plus étudiées en particulier par le groupe Max Plus ([ABG], [But]), des algèbres sur le semi-anneau tropical ([Pin], [Sim94]) et plus généralement des algèbres sur les semi-corps commutatifs idempotents ([GM02], [Zar], [Zhu]) qui  sont en fait les  <<corps>>  de caractéristique 1, comme nous le montrons ci-dessous.  Nous donnons des définitions  qui permettent de retrouver la plupart des notions introduites en algèbre tropicale tout en redonnant les notions habituelles dans le cas des corps classiques.\\
  Ce travail est motivé par l'utilisation de plus en plus diversifiée de ces algèbres qui interviennent aussi bien en Théorie du contrôle ([CGQ99], [Plus]), en Informatique théorique ([Kro], [Sim88]), qu'en Géométrie algébrique ([Izh1], [Izh2], [Mik1], [Mik2]) ou de façon plus ponctuelle dans d'autres parties des mathématiques (par exemple en Théorie des monoïdes avec des techniques de réduction en caractéristique 1 ([Cas], [Hée]) ).
 Son objectif est d'expliquer et de développer les analogies frappantes existant entre, par exemple, les résultats obtenus en algèbre Max-Plus, parfois par des méthodes d'analyse convexe, et l'algèbre linéaire classique  (voir  : [DSS], [GGB],  [CGQ04], [CGQ], [GK], [GP]...), ou encore entre les courbes tropicales et les courbes algébriques... \\
Dans ce premier article, nous revisitons  plus particulièrement l'algèbre linéaire sur les semi-corps idempotents  : le point le plus d\'elicat, 
mais le plus int\'eressant sans doute, est  la notion de point singulier vue comme g\'en\'eralisation  de celle de z\'ero. Nous retrouvons ainsi plusieurs des  <<régularités>>  déjà introduites pour les matrices et ceci nous permet de donner des définitions, nouvelles dans ce contexte, comme celle du  noyau d'une forme linéaire,  et d'utiliser la dualité pour obtenir une bonne théorie  de la dimension tropicale d'un sous-module, avec un théorème de la base incomplète et un théorème du rang.\\
Nous utiliserons aussi,  dans un  article à suivre,   ces idées  pour mettre en relief la parentée entre les courbes tropicales et les courbes alg\'ebriques habituelles (il s'agit dans les deux cas de l'ensemble des <<zéros>> d'un polyn\^ome \`a deux variables) et obtenir  de nouveaux outils algébriques pour l'étude de ces courbes.\\

{\it Remerciements.} 
   Je tiens à remercier tout particulièrement Stéphane Gaubert et Max Plus pour leur accueil, leurs encouragements et les discussions fructueuses sans lesquelles ce travail n'aurait pas été possible.

\section{ Quasi-corps}

Le but de cette section est de montrer que les semi-corps idempotents de l'alg\`ebre tropicale et les corps de l'alg\`ebre classique sont les deux facettes d'une m\^eme notion, les quasi-corps, et que l'alg\`ebre tropicale est en fait l'alg\`ebre de la "caract\'eristique 1".  

  \subsection{ Quasi-groupes}
  
  Rappelons qu'un monoïde est un ensemble muni d'une loi interne associative, admettant un \'el\'ement neutre.\\
  Un monoïde  $G$ est   {\it VN-r\'egulier } si, pour tout $x \in G$,  $ x $ appartient \`a $xGx$ (c.f. les anneaux r\'eguliers au sens de Von neumann).\\
  
    On dira qu'un \'el\'ement $x$ d'un  monoïde $(G,*)$  est {\it quasi-inversible} s'il existe un 
\'el\'ement $y$ de $G$ tel que $  x*y*x=x$ et $y*x*y=y$.  On dira alors que $x$ et $y$ sont quasi-inverses l'un de l'autre. Un {\it quasi-groupe} est un monoïde dont tous les \'el\'ements sont quasi-inversibles.\\
Il est en particulier VN-r\'egulier.\\
Réciproquement un monoïde VN-régulier est un quasi-groupe :
si $x*y*x = x$ en posant $z = y*x*y$ on a en effet $x*z*x = x$ et $z*x*z = z$.\\

 \begin{rem}
 Dans le cas commutatif $y$ est alors unique et est appel\'e le {\it quasi-opposé} de $x$ (on le notera $x^*$). \\
 En effet si on a (en notation additive), $ x+y+x=x$, $y+x+y = y$, $x+y'+x = x$, et $y'+x+y' = y'$, il en d\'ecoule : $y'+x= y'+x+y+x = y+x+y'+x = y+x$, d'o\`u $y = y+x+y =y'+x+y= y'+x+y'=y'$.\\
 \end{rem}
 
 On d\'efinit de mani\`ere naturelle les notions de {\it  sous-quasi-groupes}  et de {\it morphisme de quasi-groupes}.\\
 Deux \'el\'ements $x$ et $y$ d'un monoïde $G$ sont {\it orthogonaux}  (notation  : $x \perp y$)  si $x*z = x$ et $y*z = y$ implique $z=e$, o\`u $e$ d\'esigne l'\'el\'ement neutre de $G$. On dira que $(x_1,x_2) \in G^2$ est une {\it décomposition orthogonale} de $x \in G$ si $x = x_1*x_2$ et $x_1 \perp x_2$. Dans le cas commutatif on notera $x = x_1 \bigoplus x_2$ pour indiquer que $(x_1,x_2)$ est une décomposition orthogonale de $x$.\\

  \subsection{ Quasi-anneaux et quasi-corps}
 Un {\it quasi-anneau} est un triplet $(A,+,*)$ tel que $(A,+)$ soit un quasi-groupe commutatif, $*$ soit associative,  distributive par rapport \`a $+$, et  admette un \'el\'ement neutre, not\'e $1$ dans la suite.
 
 \begin{rem} Un quasi-anneau est un semi-anneau qui est un quasi-groupe pour l'addition.
 \end{rem}
 
 Dans la suite on notera $0$ l'\'el\'ement neutre de l'addition d'un quasi-anneau $A$, et $\epsilon$ le quasi-opposé de $1$. On aura donc, pour tout $x \in A$, $x^*=\epsilon x$.\\
 
 Si $A$ est un quasi-anneau, l'ensemble $A[X_i]_{i\in I}$ des polyn\^omes \`a coefficients dans $A$ est aussi un quasi-anneau.\\
 Un quasi-anneau est simplifiable \`a droite si $\forall (x,y,z) \in A^3$, $x*z = y*z \Longrightarrow x=y$.\\
 
 Un {\it quasi-corps} est un quasi-anneau $(K,+,*)$ tel que $(K^*,*)$ soit un groupe  (o\`u $K^*= K-\{0\})$.\\
 
\begin{ex}
Nous utiliserons plus particulièrement dans la suite les deux quasi-corps suivants : \\

Le semi-corps à deux éléments, $F_1 = \{0,1\}$, muni de l'addition telle que $0$ soit élément neutre et $1+1=1$. et de la multiplication habituelle, est un quasi-corps de caractéristique 1,  isomorphe à l'ensemble des parties d'un singleton, muni de la réunion et de l'intersection. Il est facile de vérifier que c'est le seul quasi-corps fini de caractéristique 1...

 Le semi-corps des réels max-plus, $\T$, soit sous la forme  rencontrée dans la plupart des applications, $\R \cup \{-\infty\}$ muni de la loi max comme addition et de la loi $+$ comme multiplication, soit  dans sa version, $\R_+$ muni de la loi max comme addition et de la multiplication usuelle, plus pratique pour conserver des notations algébriques générales (et plus facile à suivre par des non-spécialistes) ..

\end{ex}

  \subsection{ Caract\'eristique d'un semi-anneau}
 
 Rappelons qu'un semi-anneau $A$ se définit comme un anneau, en affaiblissant la condition $(A,+)$ est un groupe commutatif, en $(A,+) $ est un monoïde commutatif.\\
 
 Soit $A$ un semi-anneau ; on d\'efinit $H \subset \N$ comme l'ensemble des entiers $k$ tels que $k.1 + 1 = 1$. \\
 
 \begin{prop}
 Il existe un unique $ n \in \N$  tel que $H = n \N$. Cet entier  $n$ est appel\'e caract\'eristique de $A$ (notation $car(A)$).\\
 \end{prop}
  
 L'unicit\'e est \'evidente ; si $H$ n'est pas r\'eduit \`a $\{0\}$, soit $n$ le plus petit \'el\'ement non nul de $H$ : pour $m \in H$, on peut \'ecrire $ m = nk+r$ avec $0 \le r < n$ ; de $k+1=1$, on déduit par it\'eration $nk.1+1=1$ ; de $m.1 + 1=1$ on d\'eduit alors $1 = r .1+ nk.1 +1 =  r.1+1$ et donc $r=0$ par d\'efinition de $n$.\\  
  
 \begin{rem} 
 Les semi-anneaux de caract\'eristique 1 sont les semi-anneaux idempotents  (i.e. tels que $x+x=x$ pour tout $x$).
 \end{rem}
 
  On dira qu'un semi-anneau $A$ est de {\it caract\'eristique pure} p, s'il est de caractéristique $p$ et que,  pour tout $x  \not= 0$, $(k+1)x=x$ implique que $p$ divise $k$.\\
  
  \begin{prop}
  
  a) Tout semi-corps admet une caract\'eristique pure.\\
  b) Les semi-corps de caractéristique non nulle sont des quasi-corps.\\
  c) Les quasi-corps ayant une  caract\'eristique $p$ diff\'erente de $1$ sont des corps. Les quasi-corps de caract\'eristique $1$ sont les semi-corps idempotents.\\
 \end{prop}

 \noindent
 Le a) est clair.\\
 b) Si $K$ est un semi-corps et s'il existe $p \in \N$ tel que $p.1+1= 1$, soit $p.1 = 0$ (avec $ p \not=1$) et il est facile de voir que $K$ est alors un corps (de caractéristique $p$), soit $p.1 + 2.1 = 2.1$ et par itération 
 $p.1 + p.1 = p.1$, d'où, en multipliant par l'inverse de $p.1$, $1+1 = 1$ et $K$ est un semi-corps idempotent et donc bien un quasi-corps.\\  
c) Soit $K$ un quasi-corps :\\
 il existe $\alpha$ tel que $1+\alpha+1 = 1$ donc en posant $\beta = 1 + \alpha$, on a  $\beta+\beta = \beta$,  ce qui implique, puisque $K$ est un quasi-corps, $\beta = 0$, et il est alors facile de v\'erifier que $K$ est bien un corps, ou, après simplification, $1 + 1 = 1$, et $K$ est bien un semi-corps idempotent.\\  
  
  \begin{rem} M\^eme si $K$ est un quasi-corps, l'anneau de polyn\^ome $K[X]$ n'est pas simplifiable :
  si $K$ est de caract\'eristique 1,  $(X+1)(X^2+1) = (X+1)(X^2+X+1)$.\\ En particulier, il ne peut donc pas se plonger dans un quasi-corps des fractions.\\
 \end{rem}

 Plus généralement, on a la :\\
 
 \begin{prop}
 Si $A$ est un quasi-anneau de caractéristique 1,  simplifiable, pour tout couple $(x,y) \in A^2$ tel que $xy=yx$, et tout entier $n$, on a :
 $$(x+y)^n = x^n + y^n.$$
 \end{prop}
 
 En effet $(x^n+y^n)(x^n+ x^{n-1} y+ \cdots  + xy^{n-1} + y^n)=(x+y)^{2n}$.\\
 
 \subsection{Quasi-corps de caractéristique 1 et groupes ordonnés}
 
 Un quasi-anneau de caractéristique  1 est ordonné par la relation :  $a \le b$ si $a+b = b$.\\
 Un cas particulier très important est celui des quasi-corps dont l'ordre associé est total. C'est en effet le cas de tous les quasi-corps introduits en algèbre et géométrie tropicale. On parlera alors de quasi-corps totalement ordonnés.\\
Réciproquement, tout groupe totalement ordonné apparaît comme le groupe multiplicatif d'un  quasi-corps, l'addition étant donnée par $a+b = \max(a,b)$. Il suffit en fait que le groupe ait une structure de treillis.\\

   \subsection{ Modules sur un quasi-anneau}
   
   Un {\it module \`a gauche} sur un quasi-anneau $A$ est un triplet $(M,+, .)$ o\`u $(M,+)$ est un quasi-groupe, et $ . $ une loi externe de $A \times M$ dans $M$, v\'erifiant les propri\'et\'es suivantes :\\
   $\forall  a \in A$, $\forall b \in A$, $\forall m \in M$, $\forall n \in M$,  $a.(b.m) = a*b.m$, $(a+b).m = a.m+b.m$, $a.(m+n)=a.m+a.n$ et $1.m = m$.\\
  
  Si $M$ est un module  libre de base $B=(e_i)$, $x$ et $y$ appartenant à $M$ sont orthogonaux si et seulement s'ils ont des supports (relativement à $B$)  disjoints.\\

   \subsection{ Points singuliers}
   
   Il est clair que l'ensemble des applications d'un ensemble $E$ dans un quasi-groupe $G$, $F(E,G)$ est muni d'une structure de quasi-groupe, par la loi $f*g : x \longmapsto f(x)*g(x)$.\\
Soit $H$ un sous-quasi-groupe de $F(E,G)$ et $f \in H$ :  on dira que $x \in E$ est un {\it point singulier} de $f$, relativement \`a $H$, s'il existe une décomposition orthogonale de $f$ dans $H$,  $(f_1,f_2) \in H^2$,  telle que   $f_1(x)$  et $ f_2(x)$ soient  quasi-inverses l'un de l'autre.\\  En particulier si $f(x) = e$, $x$ est toujours un point singulier de $f$ ( avec pour $f_1$,  $ f$ et pour $f_2$ l'élément neutre de $H$).\\
L'ensemble des points singuliers pour $f$ relativement \`a $H$ sera not\'e $sing_H(f)$. \\
Si $H = \{0,f_1,f_2,f\}$, où $(f_1,f_2)$ est une décomposition régulière donnée de $f$, on dira alors $(f_1,f_2)$-régulière pour $H$-régulière.\\
 
Un point singulier pour un morphisme (resp : pour une application polynomiale) sera, sauf pr\'ecision contraire, un point singulier relatif à l'ensemble des morphismes consid\'er\'es (resp : des applications polynomiales).\\ L'ensemble des points singuliers d'un morphisme $f$ sera noté $Tker(f)$ (pour {\it noyau tropical} de $f$). 

Dans le cas d'un morphisme $f$ de modules ou de quasi-anneaux,  on dira que $f$ est r\'egulier s'il est r\'egulier en tout point diff\'erent de $0$, en tant que morphisme dans le quasi-groupe additif consid\'er\'e (i.e : $Tker ( f) = \{0\}$).\\ 

\begin{rem}
Dans le cas o\`u $G$ est un groupe additif, les points singuliers de $f$ sont ses z\'eros... \\
Si $f$ est un morphisme de groupes, $Tker(f)$ est son noyau.\\
\end{rem}

\subsection{Points *singuliers et zéros}

Si $f$ est un morphisme  d'un quasi-groupe $E$ dans un quasi-groupe $G$, on peut d\'efinir  une  notion duale de celle de  point singulier : on dira que $u$ appartenant \`a $E$ est un point {\it *singulier } de $f  \in F(E,G)$, ou que $f$ est *singulier en $u$,  s'il existe une d\'ecomposition orthogonale de $u$ dans $E$, $u = u_1 \bigoplus u_2$, 
telle que $f(u_1) = f(u_2)^*$.  \\
On notera $Ker^*(f)$, l'ensemble des points de $E$ *singuliers pour $f$ et on dira que $f$ est {\it *r\'eguli\`ere} si elle est *r\'eguli\`ere en tout point diff\'erent de l'\'el\'ement neutre de $E$. 

\begin{rem}Si $E$ et $G$ sont des groupes et $f$ un morphisme de groupes, les points de $E$ *singuliers pour $f$ sont encore les \'el\'ements du noyau et les deux notions co\"{\i}ncident donc.\\
  \end{rem}
  
 Dans le cas où  le quasi-groupe $G$ est additif, on dira que $u$ est un {\it zéro} d'un morphisme  $f $ de $E$ dans $G$, lorsque c'est un point *singulier de $f$.\\

 On dira de même que   $x = (x_i) \in A^{(I)}$ est un {\it zéro}  d'un polyn\^ome $P \in A[X_i]_{i \in I}$, sur un quasi-anneau commutatif $A$, si $P$ est un zéro (i.e.  un point *singulier) pour le morphisme d'évaluation en $x$,  $P \longmapsto P(x)$ (i.e. si l'on peut écrire $P =P_1\bigoplus P_2$, avec $P_1(x) = P_2(x)^*$).\\
 Cette définition s'étend sans difficulté à un point de $B^{(I)}$, où $B$ est une extension commutative de $A$, ou encore à un point d'une extension (non nécessairement commutative) dans le cas à une variable.\\
 Ceci généralise donc bien la notion habituelle de zéro d'un polynôme.\\
 En particulier $0 \in  K^n$ est un zéro  de $P \in K[X_i]$ si et seulement si le terme constant de $P$ est nul.\\
 
Les zéros d'un polynôme à une variable $P \in A[X]$, seront encore appelés {\it racines} de ce polynôme $P$.\\

Si $x \in A$ est un point singulier de l'application polynomiale $P$, $x$ est une racine de $P \in A[X]$, mais la réciproque est fausse,  deux applications polynomiales $P_1$ et $P_2$ correspondant à deux polynômes orthogonaux, n'étant pas, en général, orthogonales.\\
Cependant, sur le corps max-plus, $\T$, on peut voir facilement que ces deux notions coïncident : \\
il suffit de voir que l'inf de deux applications définies par des monômes distincts est l'application nulle ;
pour cela, on peut remarquer, si $i \not=j$, que  $a_i x^i \le a_j x^j $ pour tout $x \in \R$ implique $a_i = 0$, en faisant tendre $x$ vers $0$ ou $+ \infty$ suivant les cas...\\ 

 \section{ Alg\`ebre lin\'eaire sur les  quasi-anneaux et les quasi-corps}
 
 Dans cette section nous généralisons les principales notions de théorie des modules et des espaces vectoriels aux modules sur un quasi-corps.\\  Cependant certaines des notions habituelles admettent plusieurs généralisations distinctes...\\ 
  
  \subsection{Familles  libres, familles g\'en\'eratrices, modules libres}
  
  On définit comme habituellement les notions de familles libres, génératrices, de modules libres...
  Il faut toutefois remarquer qu'un module sur un quasi-corps  n'est pas en général libre.\\
  Par exemple sur le quasi-corps $F_1$, un module libre de base $E$ s'identifie à l'ensemble des parties finies de $E$ ([Zhu]). \\
  
  \begin{prop}
  Toutes les base d'un module libre $M$ sur un quasi-corps $K$ ont même cardinal. Ce cardinal sera encore appelé la {\it dimension} de $M$.\\
  \end{prop}
  
  Ce résultat, bien connu en caractéristique différente de $1$, se généralise facilement, puisque sur un semi-corps idempotent les éléments d'une base sont des éléments non nul ayant un support minimal. 
  
  \subsection{Groupe lin\'eaire}

  On peut montrer ([Zhu]) que le groupe linéaire sur le quasi-corps de caractéristique 1 à deux éléments, $F_1$, $Gl_n(F_1)$, est isomorphe au groupe des permutations $S_n$. \\
  Ceci provient du fait de l'unicité de la base d'un module libre (à l'ordre près) sur le quasi-corps $F_1$. Les éléments d'une base sont en effet les éléments minimaux pour la relation d'ordre. \\
  Sur un quasi-corps $K$  de caractéristique 1, les éléments de deux bases distinctes sont deux à deux colinéaires par minimalité du support... $Gl_n(K)$ sera donc formé des matrices admettant un et un seul élément non nul sur chaque ligne et chaque colonne ("$S_n \bigotimes K$").\\

  \subsection{Familles  r\'eguli\`eres et *r\'eguli\`eres}
  La notion de famille libre se révèle trop forte en général.\\ Elle correspond pour une famille $(e_i)$ à l'injectivité de l'application 
  de  $A^{(I)}$ dans $E$, \\ $(\lambda_i) \longmapsto \sum \lambda_i e_i$  ; il est naturel d'affaiblir cette condition en une condition de régularité : 
  on dira qu'une application linéaire est {\it faiblement injective } si elle est *régulière (soit  $Ker^* f = \{ 0 \})$.\\
  On dira qu'une famille $(e_i)_{i \in I}$ d'\'el\'ements d'un module $E$ sur quasi-anneau $A$ est une {\it famille r\'eguli\`ere }  (resp {faiblement libre}) si l'application de $A^{(I)}$ dans $E$, $(\lambda_i) \longmapsto \sum \lambda_i e_i$ est r\'eguli\`ere (resp {faiblement injective}).\\
  
 \begin{rem}
 La notion d'application régulière correspond, elle, à une notion de surjectivité faible... mais après transposition :  si $^t f $ est singulière en $l$, la forme linéaire  $l$ est *singulière en $f$ et donc "nulle" sur $Im f$.  On dira donc qu'un application linéaire est {\it faiblement surjective} si le sous-module orthogonal  $(Im f)'$ des formes linéaires *singulières sur $Im f$ est réduit à $\{0\}$ et qu'un famille $(f_i)$ d'élément d'un quasi-module $E$ est {\it faiblement génératice}  si le sous-module  $(\sum_i Kf_i)'$, orthogonal du sous-module engendré, est nul 
  ( On peut voir facilement que, si $E$ est un module libre, son orthogonal est réduit à $\{0\}$. Plus généralement une application linéaire entre deux modules libres est nulle si et seulement si  elle est singulière en tout point de $E$)

 Un sous-module $F$ dont l'orthogonal dans le dual est nul sera dit {\it dense}.\\
 
 Une application linéaire $f$ sera donc faiblement surjective si et seulement si  sa transposée est régulière si et seulement si son image est dense.\\

 \end{rem}

  \subsection{Matrices r\'eguli\`eres et *r\'eguli\`eres}
Soient $E$ et $F$ deux modules sur un quasi-anneau $A$. L'ensemble des morphismes (ou applications lin\'eaires) de $E$ dans $F$, $L(E,F)$ est muni d'une structure de module. 
Si $E$ et $F$ sont des modules libres de dimensions respectives $n$ et $m$, $L(E,F)$ est un module libre de dimension $mn$, isomorphe au module $M_{m,n}(A)$ des matrices \`a $m$ lignes et $n$ colonnes, \`a coefficients dans $A$.\\
On dira que  $C \in  M_{m,n}(A) $ est une {\it matrice r\'eguli\`ere } (resp {*r\'eguli\`ere}) si ses colonnes forment une famille r\'eguli\`ere (resp {*r\'eguli\`ere}) de  $M_{m,1}(A)$.\\ En appliquant les d\'efinitions, on obtient la :\\

\begin{prop}
Soit $f$ une application lin\'eaire du $A$-module libre de dimension $n$ , $E$ de base $(e_i)_{1 \le i \le n}$
dans le  $A$-module libre de dimension $m$ , de base $(f_i)_{1 \le i \le m}$ et $M$ sa matrice par rapport \`a ces deux bases. $M$ est r\'eguli\`ere (resp {*r\'eguli\`ere}) si et seulement si $f$ l'est.\\
\end{prop}

Une matrice qui n'est pas r\'eguli\`ere (resp {*r\'eguli\`ere}) sera dite {\it singuli\`ere}  (resp {*singuli\`ere}).\\
 
 \begin{prop}
 Une matrice *singulière est singulière.\\
 \end{prop}

En effet, en notant $A_i$ les vecteurs colonnes d'une matrice $A$, une relation non triviale $\sum_{I_1} x_i A_i = \sum_{I_2} x_i A_i $, où $I_1$, $I_2$ sont deux parties disjointes de $[1,n]$, donne immédiatement une décomposition $A = A_1 \bigoplus A_2$ telle que $A_1 X = A_2 X$, où $X$ est la matrice colonne des $x_i$ (en posant $x_i = 0$ pour $i \notin I_1 \cup I_2$),  $A_1$ est la matrice obtenue à partir de $A$ en remplaçant les colonnes $A_i$, $i \in I_2$ par le vecteur colonne nul, et $A_2$ la matrice complémentaire obtenue en remplaçant par le vecteur colonne nul les $A_i$ pour 
$i \notin I_2$.\\

\subsection{D\'eterminant}

Soit $K$ un quasi-corps commutatif. On d\'efinit la $K$-signature d'une permutation $\tau \in S_n$ par $k(\tau) = 1$ si la signature de $\tau $ est 1, $k(\tau) = \epsilon$ sinon.\\
Le $K$-d\'eterminant d'ordre $n$ est alors le polyn\^ome 
$$\det =  \sum_{\tau \in S_n}  k(\tau) X_{1,\tau(1)} \cdots X_{n,\tau(n)}  \in K[X_{i,j}]_{1\le i,j \le n}.$$
On dira qu'une matrice est  d-singuli\`ere lorsqu'elle est un zéro du d\'eterminant.\\ 

On pose $\det_+ = \sum_{\tau \in A_n}   X_{1,\tau(1)} \cdots X_{n,\tau(n)}  \in K[X_{i,j}]_{1\le i,j \le n}$ et\\ $\det_- =  \sum_{\tau \in S_n - An}  \epsilon X_{1,\tau(1)} \cdots X_{n,\tau(n)}  \in K[X_{i,j}]_{1\le i,j \le n}$.\\
On dira qu'une matrice $A$  est   D-singulière si $det_+(A) = det_-(A)$, c'est à dire si $A$ est $(\det_+, \det_-)$-singulière. ;  si $A$ est D-singulière elle est donc  a fortiori d-singulière, puisque c'est un zéro du déterminant..\\ 

 \begin{rem}
 Il est clair que sur un corps toutes ces notions coïncident avec la notion habituelle.\\
En algèbre tropicale le déterminant est ce qui appelé en général <<permanent>>.\\
La d-singularité  correspond  donc  en algèbre tropicale à la régularité au sens de Butkovi\v{c} (ou encore régularité tropicale ([But]).\\
De même la d-*singularité correspond à  singularité au sens de Gondran-Minoux, ou G-M-singularité (cf. [GM 84]).\\
Il est connu  que  sur le  semi-corps $\T$,  la notion de régularité au sens de Gondran-Minoux   correspond à ce que nous appelons ici  *singularité ([GB ]) et que la régularité au sens de Butkovi\v{c} correspond à la régularité tropicale, qui n'est autre que la régularité au sens ci-dessus ([Izh3], [IR])). Il est relativement aisé de voir que les démonstrations données se généralisent à tout quasi-corps totalement ordonné...
D'où la :\\

  \end{rem}

\begin{prop}
Soit  $A \in M_n(K)$ une matrice sur un quasi-corps de caractéristique 1, totalement ordonné, $K$.\\
a) $A$ est *régulière si et seulement si elle est d-*régulière.\\ 
b) $A$ est régulière si et seulement si elle est d-régulière.\\
\end{prop}

a) $A$ est *singulière signifie qu'il existe un vecteur colonne $X$ non nul et une décomposition orthogonale de $X$, $(X_1,X_2)$, telle que $AX_1 = AX_2$. Ceci correspond exactement au fait que les colonnes soient liées au sens de  G-M et il est  connu (la démonstration [GM 84 ] de se généralise sans problème) que ceci est équivalent au fait que $det_+ (A) = det_-( A)$.\\

b) Là encore il suffit de voir que la notion de matrice régulière correspond bien à celle de tropicalement régulière qui équivaut sur un quasi-corps totalement ordonné à la d-régularité en reprenant la démonstration de Z. Izhakian ([Izh3]).  \\
Pour cela on peut généraliser la construction du "revêtement" $T$ du quasi-corps $(\overline \R, max, plus)$ à un quasi-corps totalement ordonné, en posant par exemple $T = K \cup U$ où $U$ est en bijection avec $K$ par $k \longmapsto \hat{k}$ avec les règles suivantes pour $k$ et $k'$ dans $K$:\\
+ est commutative, $k + k' = max(k,k')$ si $k \not= k' $ ou $k=k'=0$, $\hat{k}$ sinon, $k + \hat{k'} = k$ si $k >k'$, $\hat{k'}$ sinon.
$\times$ est commutative, coïncide avec la multiplication sur $K$, vérifie $k \times \hat{k'} = \hat{k \times k'}$.\\
$(U \cup {0}, +, \times)$ est isomorphe, par $\hat{ }$, à $(K,+,\times)$.\\
$A \in M_n(K)$ est tropicalement singulière s'il existe un vecteur non nul $X \in K^n$ tel que $AX \in U^n$. Pour chaque ligne $i$, il existe donc au moins deux indices distincts  $j_i$ et $k_i$ tels que $a_{i,j_i} x_{j_i} = a_{i,k_i} x_{k_i}$ soient maximaux et il suffit de prendre pour $A_1$ la matrice dont tous les coefficients sont nuls sauf les $a_{i,j_i} $ et pour $A_2$ la matrice complémentaire : on a bien ainsi  $A_1 X = A_2 X$... 
Réciproquement, il est clair que si $A$ est singulier en $X$, on a bien $AX \in U$ (en calculant dans $T$).\\
 
 On peut remarquer que les matrices de déterminant non nul,  sont, en caractéristique 1, caractérisées comme suit :\\

\begin{prop}
 
 Une matrice carré $A$ a un déterminant non nul si et seulement si elle est supérieure à une matrice $S \in Gl_n$.\\

 \end{prop}

 Il suffit de remarquer que $det A$ est non nul si et seulement si l'un des produits $a_{1,\tau(1)} \cdots a_{n,\tau(n)} $ est non nul, et que ceci est vrai si et seulement si la matrice $S$ dont les seuls coefficients non nuls sont  les $a_{i,\tau(i)}$ est dans $Gl_n$.\\ 

  \section{Algèbre linéaire sur un quasi-corps de caractéristique 1 totalement ordonné}

Soit maintenant $K$ un semi-corps commutatif,  idempotent et totalement ordonné.\\

\subsection{Noyaux et Dualité}

\begin{prop}
Soit $l \in E^*$ une forme linéaire sur un module libre de dimension finie $E$.\\
$Tker l = Ker^* l$ est un sous-module de $E$. On notera plus simplement $Ker l $ ce noyau.\\
\end{prop}

Soit $B = (e_i)_{1 \le i \le n}$ une base de $E $ et $L =(a_1, \cdots, a_n)$ la matrice de $l$ dans cette base. Il est facile de voir que, le corps étant totalement ordonné, 
$x = \sum x_i e_i  $ appartient à $Tker l$ ou à $Ker^* l$ si et seulement si il existe $(i,j)$, $i \not= j$, tels que 
$a_i x_i = a_j x_j \ge a_k x_k$ pour tout $1 \le k \le n$.\\ On a donc bien l'égalité des deux noyaux dans ce cas.\\
De plus si $x$ et $y$ sont singuliers, soient $(i,j)$ et $(r,s)$ deux couples tels que $a_i x_i = a_j x_j \ge a_k x_k$  et $a_r y_r = a_s y_s \ge a_k y_k$, pour tout $1 \le k \le n$. Si les deux couples ont des supports disjoints il est clair que $a_i x_i + a_r y_r= a_j x_j  + a_s y_s \ge a_k x_k + a_k y_k$ pour tout $1 \le k \le n$, et que donc, $x+y $ est singulier.\\
Les autres cas sont encore plus simples et en en déduit facilement le résultat.\\

Pour une partie $A$ de $E$ on définit l'orthogonal $A'$ de $A$,  comme l'ensemble des formes linéaires singulières sur $A$. \\
Pour une partie $B$ de $E^*$ on définit l'orthogonal de $B$, $B^\perp$ comme l'intersection des noyaux des éléments de $B$.\\ 
On dira que $A \subset E$ est fermée si $(A')^\perp = A$.\\
 
  \subsection{Rang Tropical d'une application linéaire et Dimension Tropicale d'un sous-module}

En remarquant qu'une matrice et sa transposée ont même déterminant, ce qui précède donne immédiatement le corollaire suivant :
 
  \begin{cor}
 a) Une matrice $A \in M_n(K)$ est donc régulière (resp. *régulière) si et seulement si sa transposée l'est.\\
 b) Si $E$ est un $K$-module libre de dimension $n$, et $f \in L(E)$ est faiblement surjective, elle est aussi faiblement injective.\\
 \end{cor}

Il est naturel de généraliser la notion de famille faiblement génératrice, définie ci-dessus uniquement pour une famille faiblement génératrice d'un module libre $E$ :\\
 Dans un module libre $E$, le sous espace faiblement engendré par une famille $(f_i)$ sera l'intersection des noyaux $Ker^* l$ des formes linéaires $l \in E^*$, *-singulières sur $\sum Kf_i$.\\
Une { \it base tropicale} d'un sous-module $F$ sera alors   une famille régulière  et faiblement génératrice de $F$. On dira d'un sous-module qu'il est {\it tropicalement  libre} s'il admet une base tropicale.\\
Un sous-module tropicalement libre est donc fermé. Nous établirons la réciproque par la suite.

On a alors (voir aussi [ACG] ) :
 
 \begin{lem}
 Une matrice $A \in M_{p,n}(K)$, $p \ge n$, est régulière si et seulement si elle contient une matrice carrée extraite d'ordre $n$, régulière.\\
 \end{lem}
 
 Il est clair que si $A$ est singulière , toute les matrices carrées extraites le seront aussi.\\
 La démonstration se fait par récurrence sur $n$ : pour $n = 1$, c'est clair. En supposant le résultat vrai pour $n-1$, il existe donc une matrice carrée  d'ordre $n-1$ extraite des $n-1$ premières colonnes de $A$, régulière. Soit $(i_1, \cdots i_{n-1})$ les rangs des colonnes de cette matrice extraite et $B $ la matrice d'ordre $(n-1,n)$ extraite de $A$, dont les lignes sont celles de rangs  $(i_1, \cdots i_{n-1})$. On a donc $^t B = (L_{i_1}, \cdots , L_{i_{n-1}})$.
 En considérant la forme linéaire qui à $U = (u_1, \cdots u_n)$,  associe le déterminant de la matrice $(^t B,^t V)$ où $V = (u_{i_1}, \cdots,  u_{i_{n-1}}) $, est la matrice des composantes de $U$  correspondant aux lignes de $B$, on obtient une forme linéaire non nulle sur $K^n$,  *singulière sur les $n-1$ premières lignes de $A$. Comme $A$ est régulière, cette forme linéaire n'est pas *singulière sur toutes les lignes de $A$ et il existe donc une ligne $L_{i_n}$ telle que la matrice extraite de $A$, dont les lignes sont celles de rangs   $(i_1, \cdots i_{n})$ est régulière.\\

\begin{prop}
Soit  $E$  un $K$-module libre de dimension $n$.\\
a) Toute famille faiblement génératrice  a au moins $n$ éléments.\\
b) Toute famille régulière  a au plus $n$ éléments.\\
\end{prop}

\noindent
a) :  soit $(f_i)$ une famille faiblement génératrice d'éléments de $E$  ayant   $ k \le n-1$ éléments ; on peut supposer $k$ minimal (il est clair que si $n \ge 2$ une famille à un élément n'est pas faiblement génératrice) ; les $f_i$ sont alors non nuls ; la matrice $A$ de la famille dans une base $(e_i)$ de $E$ ne peut avoir de ligne nulle sinon il existe une forme linéaire non nulle, nulle sur tous les vecteurs de la famille. La forme linéaire $y \longmapsto det(f_1, \cdots, f_k, e_{k+1}, \cdots , e_{n-1}, y)$ est alors non nulle et singulière sur tous les $f_i$.\\
b)  : c'est immédiat par dualité.\\

\begin{rem}
  
  On peut montrer  aussi qu'une famille *régulière d'un module libre de rang $n$, a au plus $n$ éléments.\\
Ceci découle (au moins dans le cas  des réels max-plus) directement du théorème de Cramer,
et redonne cette propriété pour les familles régulières.\\
   
   \end{rem}

On peut obtenir  de plus  un <<théorème de la base tropicale incomplète>>  :\\

\begin{theo}
Toute famille régulière $(f_i)_{1\le i \le p}$  d'un module libre de rang $n$ de base $B=(e_i)$ peut se compléter en une base tropicale de $E$, en choisissant $n-p$ vecteurs dans $B$.\\
\end{theo}

Il existe une matrice extraite régulière d'ordre $k$ de $A$,  la matrice dans la base $B$ de la famille $(f_i)$. En complétant par les $e_i$ correspondant aux lignes n'apparaissant pas dans la matrice extraite, on obtient une famille régulière.\\

 De tout ceci, on déduit finalement la :

\begin{prop}
Soit $F$ un sous-module tropicalement  libre d'un module libre $E$ de dimension  finie $n$. Deux bases tropicales  de $F$ ont même cardinal. Cet entier est la {\it dimension tropicale} de $F$.
\end{prop}

Soit  maintenant  $(f_i)_{1 \le i \le k}$ une base tropicale de $F$, $(e_i)_{1\le i\le n} $ une base de $E$ et $A$ la matrice des composantes des $f_i$ dans la base $(e_i)$. D'après la proposition précédente $k \le n$ et si $k = n$, la famille est tropicalement génératrice et donc $F = E$. On peut donc supposer $k < n$.\\
  On peut donc choisir des vecteurs $f_{k+1}, \cdots, f_n$ dans la base $B$, tels que la famille $(f_i)$ soit une base tropicale de $E$. Si $G$ est le sous-module engendré par les vecteurs $f_{k+1}, \cdots, f_n$, $G$ est un module libre et $F \bigoplus G$ est dense dans $E$.\\
  Si $(g_i)_{1 \le i \le k'}$ est une autre base tropicale de $F$,  $(g_i)_{1 \le i \le k'} \cup (f_{k+1}, \cdots, f_n)$  est une famille génératrice et tropicalement libre de $E$ ce qui implique $k = k'$.\\

 \begin{rem}
 On peut définir plus généralement  la {\it dimension  tropicale} d'un sous-module $F$ d'un module libre $E$ de rang fini, comme le cardinal maximal d'une famille régulière d'éléments de $F$ (notation : $dim_T(F)$). \\  
 Pour un module libre la dimension tropicale est bien entendu égale à la dimension.
De plus cette définition de la dimension tropicale   est cohérente avec  la définition  de rang tropical d'une matrice donnée  par Z.  Izhakian :
en effet si $A$ est une matrice, le {\it rang tropical} de $A$ est défini  comme l'ordre maximal d'une matrice carrée régulière extraite de $A$. Or ceci est bien la dimension tropicale de l'image de $A$.\\
 \end{rem}

\begin{lem}
Si $(f_i)$ est une famille régulière maximale dans un sous-module $F$ d'un module libre $E$ de dimension $n$. Pour tout $g \in E$ tel que la famille $(f_i) \cup (g)$ soit régulière, il existe une forme linéaire *singulière sur les $f_i$ et pas sur $g$. Le sous-module $H$ tropicalement engendré par la famille $(f_i)$ contient donc $F$.\\
\end{lem}

Si $g$  est tel que la famille $(f_i) \cup (g)$ soit régulière, on peut compléter cette famille en une base tropicale de $E$ par des $f_j$ pris dans une base $B$, et la forme linéaire \\ $ y \longmapsto det(f_1, \cdots, f_k, y, f_{k+1}, \cdots, f_{n-1})$ serait *singulière sur les $f_i$  et non *singulière en $g$ puisque la famille est régulière.\\

On en déduit que : 

\begin{prop}
a) La dimension tropicale d'un sous-module $F$ est aussi celle du module tropicalement  engendré.\\
 b) Les sous-modules ayant une famille finie tropicalement génératrice sont les sous-modules  tropicalement libres : toute famille régulière maximale est une base tropicale.\\
 c) Plus généralement un sous-module  $F$ est tropicalement libre si et seulement si il est fermé.\\
\end{prop}

\noindent
a) Si $(f_i)$ est une famille régulière maximale dans $F$, par le théorème de la base incomplète, on peut choisir une partie de la base $B$ de $E$, complétant la famille $(f_i)$ en une base tropicale de $E$. Le sous-module libre $G$ ainsi construit est de rang $n-k$, où $k$ est le cardinal de la famille $(f_i)$ ; il est en somme directe avec $F$ ( si $f \in F \cap G$ était non nul, la famille $(f_i) \cup (f) $ serait régulière...). Il est alors, d'après le lemme précédent, en somme directe avec le sous module $H$ tropicalement engendré par $F$. On en déduit immédiatement que le rang tropical de $H$ est encore $k$.\\ 
b) et c)  Si $F = H$,  $(f_i)$ est donc une base tropicale de $F$.\\

On obtient une caractérisation des sous-modules tropicalement libres :

\begin{prop}
 a) Les sous-modules tropicalement libres sont les intersections finies de noyaux de formes linéaires.\\
b) Si $f$ est une application linéaire entre deux modules libres de dimensions finies $E$ et $F$, $Tker f$ est un sous-module tropicalement libre de $E$.\\
 \end{prop}

a) Par dualité il est clair que $F'$, l'ensemble des formes linéaires *singulières sur $F$, est tropicalement  libre et si $F$ est tropicalement libre, il  est alors égal à l'intersection des noyaux des éléments d'une base tropicale de $F'$.\\
Réciproquement, l'intersection des noyaux d'une famille $(l_i)$ du dual de $E$ est l'orthogonal du sous-espace engendré et est donc tropicalement libre.\\
b) $Tker f$ est l'intersection des noyaux des formes linéaires $p_i \circ f$, où les $p_i$ sont les formes coordonnées d'une base de $F$.\\

Pour les applications linéaires, ceci donne un théorème du rang :

\begin{theo}
Soit $f$  une application linéaire entre deux $K$-modules libres, de dimensions finies, $E$ et $F$.\\
  On a :
  $$\dim E =  rg_T(f) + \dim_T(Tker f)$$ où $rg_T(f) = \dim_T(Im f)$ est le rang tropical de $f$.\\
\end{theo}

 Pour ceci on considère un <<supplémentaire tropical>>  $G$ de $Tker f$ (obtenu en complétant une famille régulière maximale de $Tker f$ en une base tropicale de $E$ , avec des vecteurs pris dans une base de $E$) . La restriction $g$ de $f$ à $G$, est alors une application du module libre $G$ dans $F$ d'image $Im f$.
 De plus elle est clairement régulière et on a donc bien $\dim G = \dim_T Im f$. \\

      \begin{center}
{\large  REFERENCES}
\end{center}

 [ABG] M. AKIAN, R. BAPAT, S. GAUBERT.{\it Max-plus algebras},  Handbook of Linear Algebra (Discrete Mathematics and Its Applications, L.HOGBEN ed.), Chapter  25, vol. 39, Chapman  \& Hall / CRC, 2006. \\ 
 
 [AGG] M. AKIAN, S. GAUBERT, A. GUTERMAN. Matrix ranks and linear independence over tropical semi-ring. Preprint.\\

 [BG] L.B. BEASLEY, A.E. GUTERMAN. Rank inequalities over semi-rings. J. Of Korean Math. Soc. 42 (2), 2005, 223-241.\\
 
  [BP] L.B. BEASLEY, N. J. PULLMAN,  Semi-ring ranf versus column rank. Linear Algebra Appl. 101 (1988) 33-48.\\

 [But] P. BUTKOVIC. Max-algebra: the linear algebra of combinatorics?. Linear Algebra and Appl., 367, 
2003, p. 313-335. \\

[Cas] A. CASTELLA. Lawrence-Krammer-Paris representation under graph automorphisms. ArXiv:math 0803.1115 (2008).\\

[CGQ99] G. COHEN, S. GAUBERT, J.-P. QUADRAT. Max-plus algebra and system theory : where we are and where to 
go now. Annual Reviews in Control, 23 , 1999, p. 207–219. \\

[CGQ04] G. COHEN, S. GAUBERT, J.-P. QUADRAT.  Duality and separation theorems in Idempotent semimodules.
"Linear Algebra and Appl.,  379,  2004,  p. 395–422.\\

[CGQ] G. COHEN, S. GAUBERT, J.-P. QUADRAT. Regular matrices in max-plus algebra. Preprint.\\

[CGB] R.A. CUNINGHAME-GREEN, P. BUTKOVIC. Bases in max-algebra.  Linear Algebra Appl. 389 (2004) 107-120.\\

[DSS] M. DEVELIN, F. SANTOS, B. STURFELS. On the rank of a tropical matrix. Discrete and Computational Geometry (E. Goodman, J. Pach, and E. Welzl eds.) MRSI Publications, Cambridge Univ. Press, 2005.\\

[GB] S. GAUBERT, P. BUTKOVIC. Sign-nonsingular matrices with unbalanced determinant in symetrised semirings.  Linear Algebra Appl. 301 (1999) 195-201.\\.\\

[GG] S. GAUBERT, J. GUNAWARDENA. The duality theorem for min-max functions. C.R.A.S.  Paris 326, Série I (1998) 43-38.\\

 [GK]    S. GAUBERT, R. KATZ. Max-Plus Convex Geometry. Preprint.\\
 
 [GP] S. GAUBERT, M. PLUS. Methods and applications of (max,+) linear algebra. STACS'97, LCNS 1200 (1997).\\

[GM77]  M. GONDRAN, M. MINOUX. Valeurs propres et vecteurs propres dans les dioïdes et leur interprétation en 
théorie des graphes. EDF, Bulletin de la Direction des Etudes et Recherches, Serie C, Mathématiques 
Informatique",  2, 1977, p. 25-41. \\

[GM84]  M. GONDRAN, M. MINOUX. Linear algebra in dioids : a survey of recent results. Annals of Discrete 
Mathematics, 19, 1984, p. 147-164. \\

[GM02]  M. GONDRAN, M. MINOUX. Graphes, Dioïdes et semi-anneaux, TEC \& DOC, Paris, 2002. \\

[Hée] J.Y. H\'EE. Une démonstration simple de la fidélité de la représentation de Lawrence-Krammer-Paris. Preprint.\\

[Izh1] Z. IZHAKIAN. Tropical arithmetic and algebra of tropical matrices. ArXiv:math. AG/0505458, 2005.\\

[Izh2] Z. IZHAKIAN. Tropical varieties, ideals and an algebraic nullstensatz. ArXiv:math. AC/0511059, 2005\\

[Izh3] Z. IZHAKIAN. The tropical rank of a tropical matrix. ArXiv:math. AC/0604208, 2006.\\

[IR] Z. IZHAKIAN, L. ROWEN. Supertropical algebra. ArXiv:math AC/0806.1171 (2008).\\

[KR] K.H. KIM, F.W. ROUSH. Kapranov rank versus tropical rank. ArXiv:math. CO/0503044 v2.\\

[Kro] D. KROB The equalities problem for rational series with multiplicities in the tropical semi-ring is undecidable. Proceedings of ICALP 92 (W. Kuich ed) LNCS 623 101-112 (1992).\\

[Mik1]  G. MIKHALKIN. Amoebas of algebraic varieties and tropicalgeometry. Different faces of geometry, 
Int. Math. Ser.(N.Y.), vol.3, Kluwer/Plenum, NewYork,  2004, p. 257–300.\\

[Mik2]  G. MIKHALKIN. Enumerative tropical algebraic geometry in $\R^2$.  J. Amer. Math. Soc. 18 (2005) 313-377.\\

[Pin] J.-E. PIN. Tropical Semirings. Idempotency (J. Gunawerdana ed.), Publications of the Isaac 
Newton Institute, Cambridge University Press, 1998. \\

[Plus] M. PLUS. Linear systems in (max, +)-algebra. Proceedings of the 29th Conference on Decision and Control, Honolulu, Dec. 1990
p. 177-207. 
. \\

[RST] J. RICHTER-GEBERT, B. STURMFELS, T. THEOBLAND. First steps in tropical geometry. Contemporary Mathematics, 377, Amer. Math. Soc. (2005), 289-317.\\

[Sim78] I. SIMON. Limited subsets of the free monoid. Proc. of the 19th Annual Symposium on Foundations of 
Computer Science,  IEEE, 1978, p. 143–150. \\

[Sim88] I. SIMON.  Reconizable sets with multiplicities in the tropical semi-ring. Proceedings of MFCS 88, (M.P. Chytil et al. eds) LNCS, 324, 107-120..\\

[Sim94] I. SIMON. On semigroups of matrices over the tropical semiring. Theor. Infor. and Appl.,  28
 (1994), p. 277–294. \\

[Zar] K. ZARETSKI.  Regular elements in the semi-group of binary relations. Uspeki Mat. Nauk 17(3) (1962) 105-108.\\

[Zhu]  Combinatorics and characteristic one algebra, Preprint (2000).\\

  \end{document}